\documentclass{gtart}     
\usepackage{amssymb,graphicx} 
\usepackage{amsmath,amscd} 

\usepackage{psfrag}
\def\psfraga <#1,#2> #3#4{\psfrag {#3}{\rlap{\kern #1 \raise #2\hbox{#4}}}}

\input rlepsf
\let\relabela\adjustrelabel

\input gtmonout
\volumenumber{2}
\volumeyear{1999}
\volumename{Proceedings of the Kirbyfest}
\pagenumbers{201}{213}
\papernumber{11}
\received{22 March 1999}\revised{23 August 1999}
\published{18 November 1999}

\newtheorem{thm}{Theorem}
\newtheorem{lem}[thm]{Lemma}         
\newtheorem{claim}[thm]{Claim}     
\theoremstyle{definition}
\newtheorem*{rem}{Remark}           
\newtheorem{example}[thm]{Example} 
\begin{document}

\title{Configurations of curves and geodesics on surfaces}
\authors{Joel Hass\\Peter Scott}
\address{Department of Mathematics, University of California\\Davis, CA 95616,
USA}
\secondaddress{Department of Mathematics, 
University of Michigan\\Ann Arbor, MI 48109, USA}
\email{hass@math.ucdavis.edu, pscott@math.lsa.umich.edu}                  

\begin{abstract}
We study configurations of immersed curves in surfaces and surfaces in
3--manifolds. Among other results, we show that primitive curves have only
finitely many configurations which minimize the number of double points. We
give examples of minimal configurations not realized by geodesics in any
hyperbolic metric.
\end{abstract}
\makeshorttitle


Let $f$ and $g$ be general position immersions of a manifold $M$ into
the interior of a manifold $N$. We will say that $f$ and $g$ {\em have
the same configuration} if there is a regular homotopy from $f$ to $g$
through general position immersions. Equivalently, there is an ambient
isotopy of $N$ moving $f(M)$ to $g(M)$. This defines an equivalence
relation on general position immersions, and an equivalence class will
be called a {\em configuration}. All the immersions in a configuration
``look the same'' in a precise sense. In this paper we will be
interested in the cases when the dimension of $M$ is $1$ or $2$ and
the dimension of $N$ is $2$ or $3$. Our aim is to explore the question
of how many configurations a given homotopy class can have. For
primitive curves on a surface, we show that the number is finite if
one restricts to immersions with the least possible number of double
points, but little can be said for curves with excess intersections.
We then consider the possible configurations of closed geodesics on a
surface equipped with a hyperbolic metric. It is well known that
geodesics in a hyperbolic metric minimize the number of double points
in their homotopy class. It was shown by Shephard \cite{S} and
Neumann-Coto \cite{N} that any curve configuration (not necessarily
connected) which minimizes the number of double points is realized by
shortest geodesics in some metric. We construct examples which show
that some configurations cannot be realized by closed geodesics in a
hyperbolic metric.

We say that a map of the circle into a surface is {\em primitive} if it is
not homotopic to a proper power of some other map. Our first and most
general result is the following.

\begin{thm}
\label{finite} Let $f\co S^{1}\rightarrow F$ be a primitive map of the circle
into an orientable surface. Then the general position immersions which are
homotopic to $f$ and have the minimal possible number of double points
belong to only finitely many configurations.
\end{thm}

\proof If $f$ is nullhomotopic in $F$, then any general position
immersion which is homotopic to $f$ and has the minimal possible number of
double points must be an embedding. Further such an embedding will bound a
disc. It follows that there are two configurations possible for such maps,
one for each orientation of the curve. Thus the theorem holds for the
2--sphere. If $f$ is homotopically essential in $F$, and is homotopic to an
embedding, there is only one configuration possible among embeddings
homotopic to $f$. For the torus $T$, the assumption that $f\co S^{1}\rightarrow
T$ is primitive implies that $f$ is homotopic to an embedding and so has a
unique configuration.

Assume now that $\chi (F)<0$ and pick a hyperbolic metric on $F$. The
pre-image of $f(S^{1})$ in the universal cover $\mathbb{H}^{2}$ consists of a
line (in the topological sense) $l$ and its translates $\{gl\}$, $g\in G$.
These lines will not be geodesics in general, but each will lie in a bounded
neighborhood of a unique geodesic. As $f$ represents a primitive element of $%
\pi _{1}(F)$, no two of these lines have the same endpoints. The minimality
of the number of double points of $f$ implies that any two of these lines
meet in at most one point, as is the case with hyperbolic geodesics. Let $%
p_{ij}$ denote the point of intersection of distinct translates $l_{i}$ and $%
l_{j}$, with the convention that $p_{ij}$ does not exist if $l_{i}$ and $%
l_{j}$ are disjoint.

\begin{claim}\label{claim1}
If we know the side of $l$ on which $p_{ij}$ lies for all $i,j$, then the
configuration of $f$ is determined.
\end{claim}

\proof Note that the assumption in the claim implies that for each $%
l_{k}$ we know the side of $l_{k}$ on which $p_{ij}$ lies for all $i,j$. We
will construct the configuration of lines one at a time, starting with $%
l=l_{1}$. Assume that the lines $l_{1,}\ldots ,l_{n-1}$ have a unique
configuration. We will establish that the configuration of the lines $%
l_{1,}\ldots ,l_{n}$ is also unique. Consider the choices when we add the
additional line $l_{n}$. Two disjoint lines in $\mathbb{H}^{2}$ cannot be
interchanged by a homotopy of $f$, unless they have the same endpoints, as
they do not lie within a bounded distance of one another. But the assumption
that the curve $f$ is primitive implies that no two lines have the same
endpoints. Hence if $l_{n}$ is disjoint from $l_{i}$ then the side of $l_{i}$
on which it lies is determined. Suppose that $l_{n}$ crosses some $l_{k}$.
The points $p_{ik},i<n$, in which the previous lines meet $l_{k}$, divide $%
l_{k}$ into several arcs. As we know on which side of $l_{i}$ the point $%
p_{kn}$ lies, we know in which of these arcs $p_{kn}$ lies. It follows that
up to an isotopy of the lines $l_{1,}\ldots ,l_{n}$, there is at most one
possible way in which to add $l_{n}$. Now induction on $n$ shows that the
collection of all translates of $l$ is determined up to ambient isotopy of $%
\mathbb{H}^{2}$. Further, if we have two immersions $f$ and $g$ of $S^{1}$ in $%
F $ such that the corresponding families of lines in $\mathbb{H}^{2}$ are
ambient isotopic, we claim that the isotopy can be chosen to be equivariant
under the action of $\pi _{1}(F)$ on $\mathbb{H}^{2}$, so that $f$ and $g$ must
have the same configuration as claimed. The way to do this is first to
ensure that the isotopy is equivariant when restricted to the intersection
points of the two families of lines, then to ensure equivariance of the
isotopy when restricted to the union of the lines and finally to ensure that
the entire isotopy is equivariant by defining it equivariantly on each of
the regions into which the union of the lines divides the hyperbolic
plane.\endproof

\begin{claim}
Let $\gamma $ be a closed geodesic in some hyperbolic metric on $F$, so that 
$l$ and its translates are geodesics in the hyperbolic plane $\mathbb{H}^{2}$.
Fix a line $l_{i}$ which crosses $l$. Then the number of lines which cross
both $l$ and $l_{i}$ is finite.
\end{claim}

\proof The entire configuration of lines projects to a closed curve
in $F$, which must have only finitely many double points, and it follows
that there are only finitely many values for the angles between any two
lines $l_i$ and $l_j$ which meet. In particular, the angles are bounded
uniformly away from zero. This yields an upper bound to the lengths of the
sides of any triangle formed by these lines. If the number of $l_{j}$'s
which cross both $l$ and $l_{i}$ is not bounded, then since the set of all
lines cannot accumulate, there must be triangles of unbounded size, a
contradiction.\endproof

Now we can complete the proof of Theorem \ref{finite}. We return to
the general situation where no metric is assumed. As any two of the
lines forming the pre-image of $f$ intersect in at most one point, the
intersections of these lines correspond to the intersections of the
corresponding geodesics in $\mathbb{H}^{2}$. It follows that the
conclusion of the preceding claim applies, so that the number of
$l_{j}$'s which cross both $l$ and $l_{i}$ is a finite number
$m_{i}$. Now the translates of $l$ which cross $l$ fall into $2n$
orbits under the action of the stabiliser of $l$, where $n$ denotes
the number of double points of $f$. Let $l_{1},...,l_{2n}$ denote one
representative from each orbit. The total number of points
$p_{ij},1\leq i\leq 2n$, such that $l_{j}$ crosses $l$ and $l_{i}$ is
$ m=m_{1}+m_{2}+...+m_{2n}$. For each $l_{j}$ which crosses both of
$l$ and $l_{i}$, there are at most two choices for which side of $l$
the point $p_{ij} $ occurs. Hence the total number of choices for
which side of $l$ these $m$ points lie is bounded by $2^{m}$. But
these choices determine completely on which side of $l$ every $p_{ij}$
lies. Now Claim \ref{claim1} implies that there are at most $2^{m}$
possible configurations. This completes the proof of Theorem
\ref{finite}.\endproof

\begin{rem}
The condition that the number of double points be minimal is essential for
Theorem \ref{finite}.
\end{rem}

Even if one restricts the number of double points to two, there is a curve
on a surface whose homotopy class contains infinitely many distinct
configurations. An example, as shown in Figure 1, can be obtained by
beginning with a simple closed curve $C$ on a surface $F$, choosing a simple
arc $\lambda $ on $F$ which meets $C$ only in its endpoints, and isotoping a
small arc of $C$ at one end of $\lambda $ until it runs back and forth along 
$\lambda $ and cuts $C$ twice near the other end of $\lambda $. For most
surfaces $F$, the relative homotopy class of $\lambda $ can be chosen in
infinitely many different ways, yielding infinitely many distinct
configurations with two double points which are all homotopic to the initial
embedding.

\begin{figure}[ht!]
\centering
\includegraphics[width=.5\textwidth]{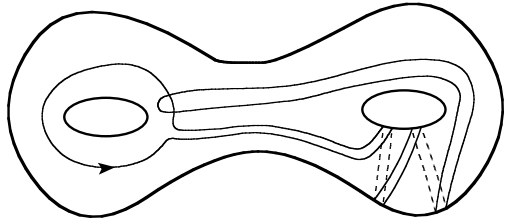}
\caption{Extra double points on a surface of genus two}
\label{finger}
\end{figure}

We will now consider some examples which show that configurations of curves
on a surface with minimal self-intersection cannot always be realized by a
geodesic in a hyperbolic metric.

\begin{lem}\label{lemma1}
Let $g_{t}$ be a family of Riemannian metrics on a closed manifold, let $%
\gamma $ be a closed curve in $M$, and let $\gamma _{t}$ be a shortest
closed geodesic homotopic to $\gamma $ in the metric $g_{t}$. If $\gamma
_{0} $ is the unique geodesic in its homotopy class then $\gamma _{t}$
varies continuously with $t$ at $t=0$. If each $\gamma _{t}$ is unique, then
the whole family is continuous.
\end{lem}

\proof Let $N_{\epsilon }$ be an $\epsilon $--neighborhood of $%
\gamma _{0}$ in the metric $g_{0}$. If there are $\gamma _{t_{i}}$ not
entirely contained in $N_{\epsilon }$ for a sequence $t_{i}\rightarrow 0$,
then a subsequence of these converges by an application of Ascoli's Theorem,
and the limit will be a geodesic not entirely contained in $N_{\epsilon }$,
but homotopic to $\gamma _{0}$ and having the same length. Thus $\gamma _{t}$
lies inside $N_{\epsilon }$ for $t$ sufficiently small, and the family of
geodesics varies continuously at $t=0$. \endproof

Our first example, for simplicity of construction, considers intersections
of three simple curves. We then describe a similar, but more complicated
example which uses a single singular curve.

\begin{example}
There are three simple closed curves on a punctured torus $F$ which have
several minimal intersection configurations, of which only one is achieved
by geodesics in any hyperbolic metric on $F$.
\end{example}

\begin{figure}[ht!]
\centering
\begin{minipage}{.45 \textwidth}
\psfraga <-2pt,-3pt> {b}{$\beta$}
\psfraga <-2pt,0pt> {a}{$\alpha$}
\psfraga <-2pt,-2pt> {g}{$\gamma$}
\centerline{\small
\includegraphics[width=.9\textwidth]{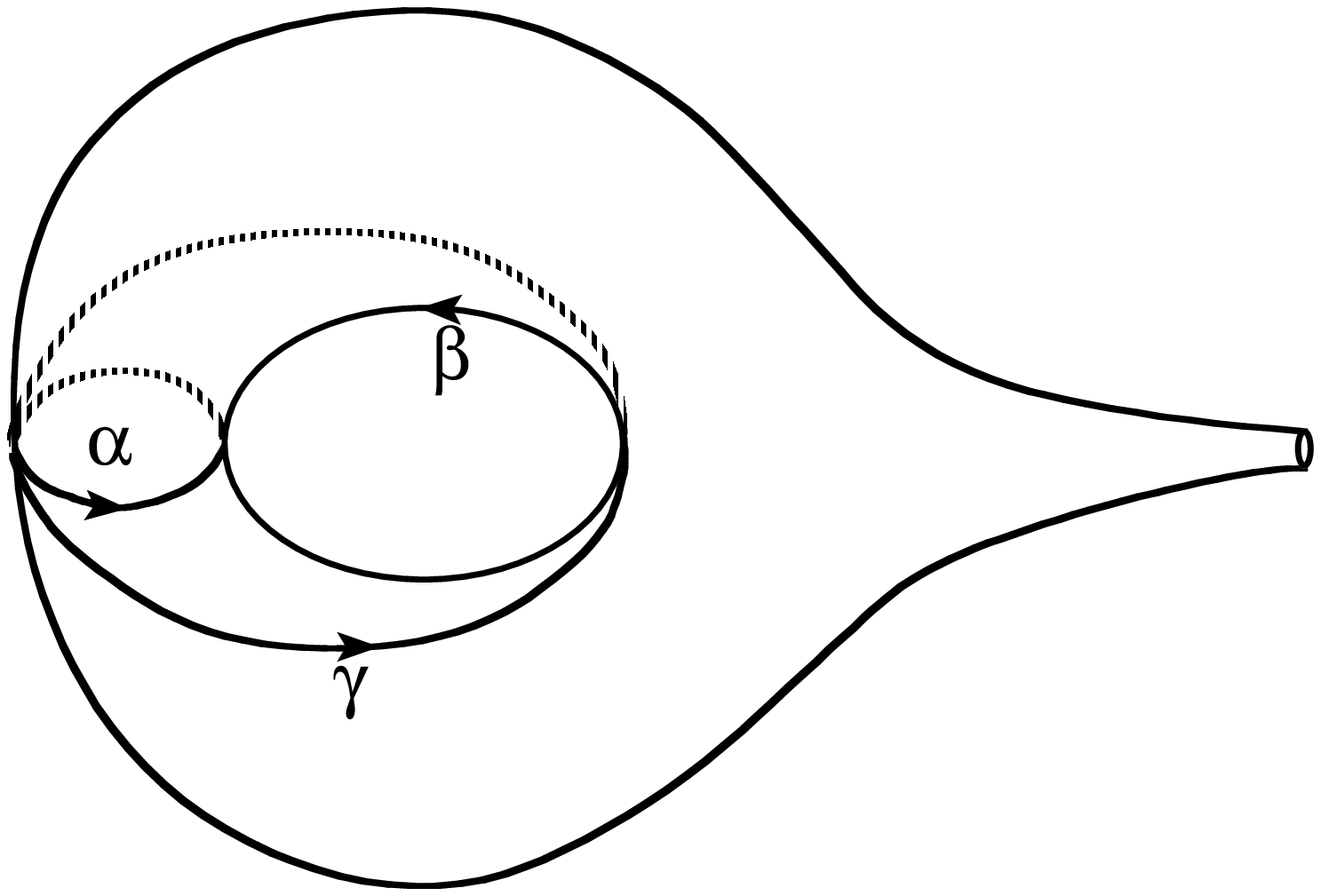}}
\end{minipage}
\begin{minipage}{.45 \textwidth}
\psfraga <-2pt,-4pt> {b}{$\beta$}
\psfraga <-4pt,0pt> {a}{$\alpha$}
\psfraga <-2pt,0pt> {g}{$\gamma$}
\centerline{\small
\includegraphics[width=.9\textwidth]{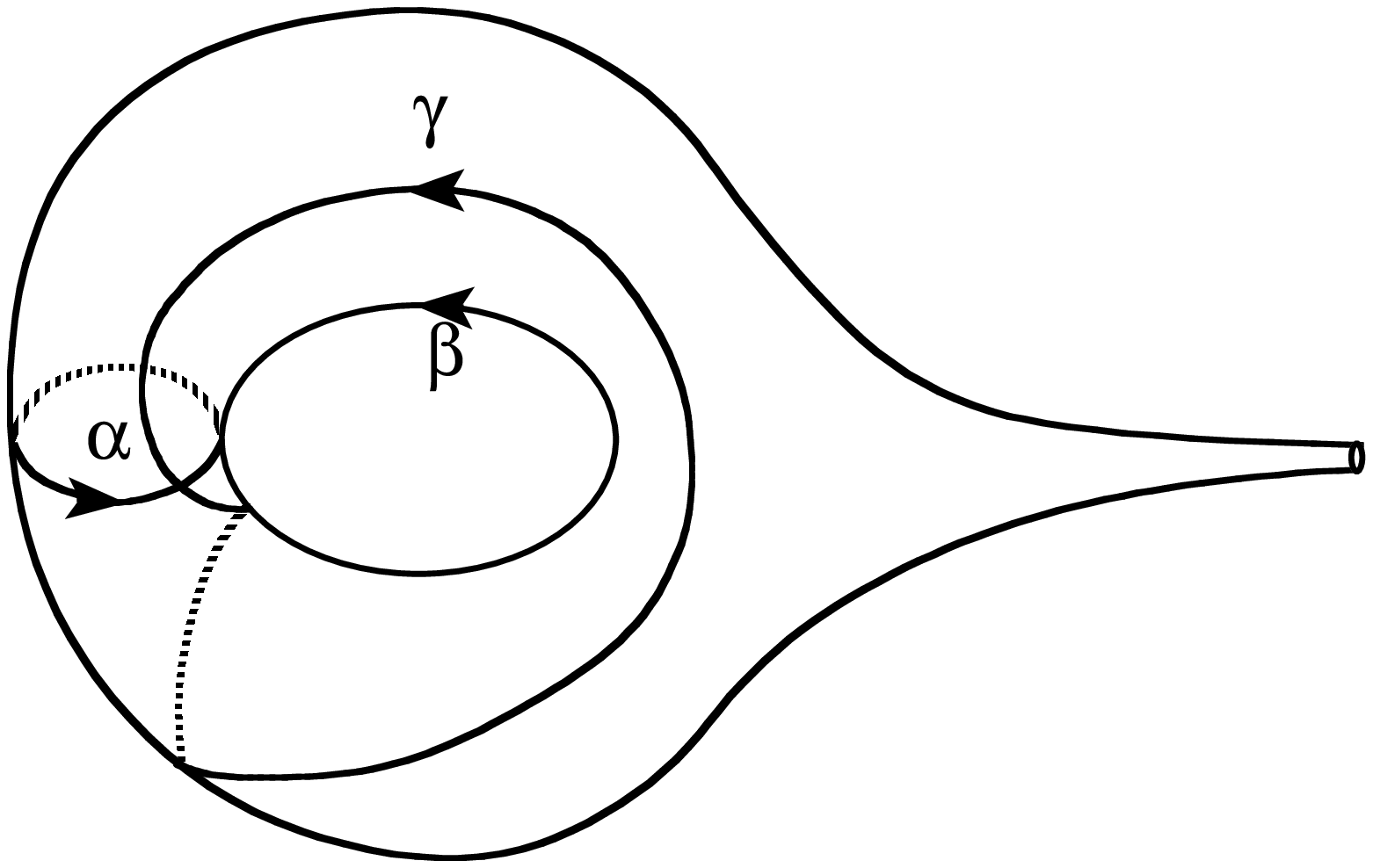}}
\end{minipage}
\caption{Only the first configuration can be realized by hyperbolic
geodesics.}
\label{3curves}
\end{figure}

Let $a$ and $b$ be a basis for $\pi _{1}(F)$ representing a longitude and
meridian, and let $\alpha $, $\beta $ and $\gamma $ be closed geodesics
representing $a$, $b$ and $ab$. Each of these curves is simple and each pair
cross in a single point. See Figure~\ref{3curves}.

The punctured torus has an involution $\tau \co F\rightarrow F$ which fixes
three points and such that $\tau (a)=a^{-1},\tau (b)=b^{-1}$. So $\alpha $
and $\beta $ are preserved by the involution. We have 
\[
\tau (ab)=a^{-1}b^{-1}=a^{-1}(b^{-1}a^{-1})a=a^{-1}(ab)^{-1}a 
\]
so that $ab$ is taken to a conjugate of its inverse, and the geodesic $%
\gamma $ is also preserved. Hence each curve is invariant, but reversed, and
so its image contains two fixed points of $\tau $. For any pair of the three
curves, the unique point at which they intersect must be fixed by the
involution.

If all three of $\alpha $, $\beta $ and $\gamma $ intersect at a
common point $x$ then this point is fixed by $\tau $, as are three
additional and distinct points, one on each of $\alpha ,\beta $ and
$\gamma $. This would result in more than three fixed points for $\tau
$, a contradiction. So $\alpha \cap \beta ,\beta \cap \gamma $ and
$\alpha \cap \gamma $ are three distinct points on $F$, as in
Figure~\ref{3curves}. Now suppose that there is more than one possible
configuration in the homotopy class of the three curves, realized by
two distinct hyperbolic structures $T_{1}$ and $T_{2}$.  We can
connect the two structures in Teichmuller space by a path of
hyperbolic structures $T_{t}$. By Lemma \ref{lemma1}, closed geodesics
in a given homotopy class vary continuously on the surface as we
follow a path of hyperbolic metrics in Teichmuller space. The above
argument shows that for each metric $T_{t}$, the unique geodesics in
the homotopy classes $a,b$ and $ab$ have no triple points. It
follows that we cannot change configurations.  However there is a
complementary region of these three curves which is a triangle---in
fact two of them are. So topologically it is possible to alter the
configuration by sliding one of the edges of this triangle across the
opposite vertex. The resulting configuration still minimizes the
number of intersection points but cannot be realized in any hyperbolic
metric. See Figure~\ref{3curves}. The same example can be put into any
surface, by constructing it inside a subsurface homeomorphic to a
torus with a disk removed.

\begin{example}
A connected closed curve on a punctured torus with several minimal
intersection configurations, of which only one is achieved by a geodesic in
a hyperbolic metric.
\end{example}

\begin{figure}[ht!]
\centering
\begin{minipage}{.4 \textwidth}
\centerline{
\relabelbox\small
\epsfxsize 1.2\textwidth\epsfbox{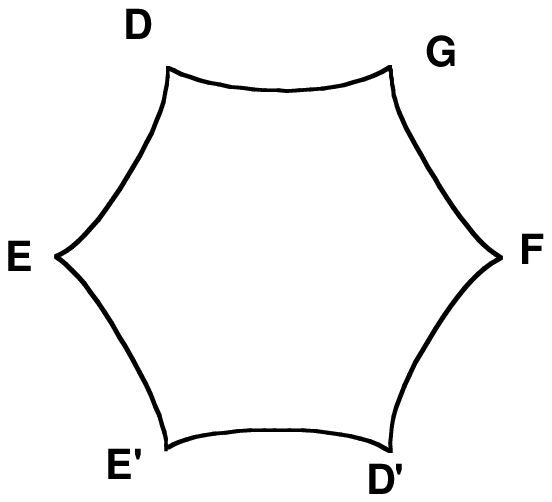}
\relabela <-2pt,-2pt> {D}{$D$}
\relabela <-2pt,-2pt> {D'}{$D'$}
\relabela <-2pt,0pt> {E}{$E$}
\relabela <-2pt,-3pt> {E'}{$E'$}
\relabela <-3pt,0pt> {G}{$G$}
\relabela <-2pt,0pt> {F}{$F$}
\hglue -4cm}
\endrelabelbox
\end{minipage}
\begin{minipage}{.59 \textwidth}
\centerline{
\relabelbox\small
\epsfxsize .9\textwidth\epsfbox{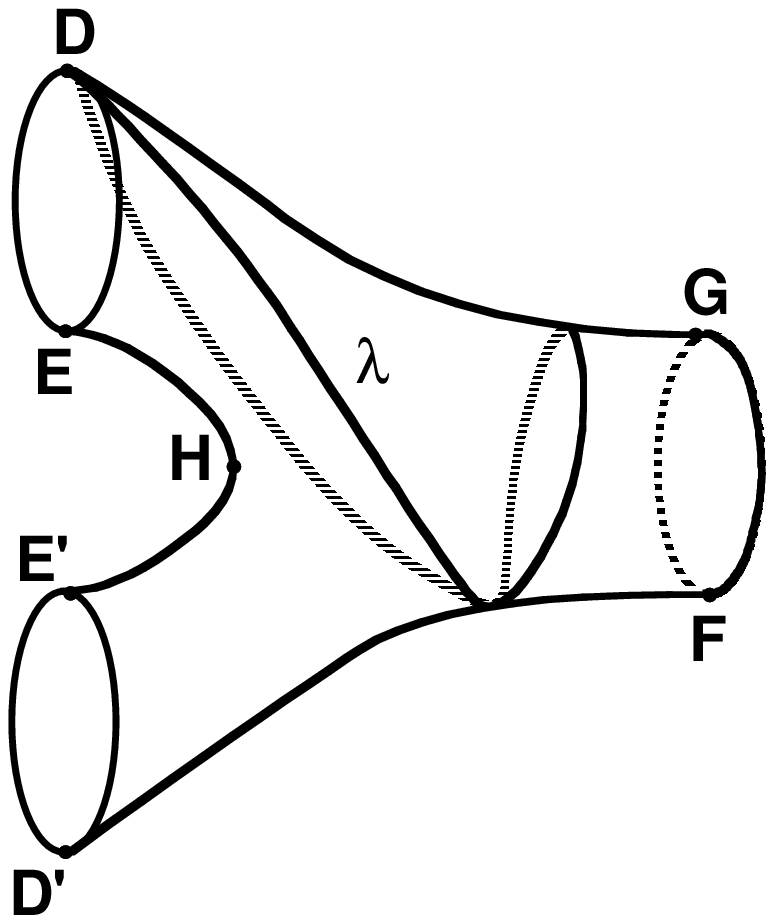}
\relabela <-2pt,0pt> {D}{$D$}
\relabela <-2pt,0pt> {D'}{$D'$}
\relabela <-2pt,0pt> {E}{$E$}
\relabela <-2pt,0pt> {E'}{$E'$}
\relabela <-2pt,0pt> {G}{$G$}
\relabela <-2pt,0pt> {F}{$F$}
\relabela <-2pt,0pt> {H}{$H$}
\relabela <-2pt,0pt> {l}{$\lambda$}
\endrelabelbox
\hglue -1cm}
\end{minipage}
\caption{An all right hexagon and a geodesic arc in a pair of pants}
\label{1curve}
\end{figure}

We start with an all right angled hyperbolic hexagon $DEE^{\prime }D^{\prime
}FG$, then double it along the edges $EE^{\prime },D^{\prime }F$ and $GD$ to
obtain a pair of pants $X$ with a hyperbolic metric, as in Figure~\ref{1curve}.
Thus $X$ admits a reflection involution $\sigma $ which
interchanges the two hexagons. It also admits an orientation preserving
involution $\tau $ which fixes a single point $H$ of $X$, where $H$ is the
midpoint of the arc $EE^{\prime }$. Now choose a geodesic loop $\lambda $ on 
$X$ based at $D$ as shown in Figure~\ref{1curve}. This loop is not a closed
geodesic, as there will be a corner at $D$. It is freely homotopic to the
square of the boundary component which contains $F$ and $G$, so it cannot be
simple. However, it can be realized with only one double point and hence has
exactly one double point. As each boundary component of $X$ is preserved by $%
\sigma $ but reversed in orientation and as $D$ is fixed by $\sigma $, it
follows that $\lambda $ is preserved by $\sigma $ but with reversed
orientation. Hence $\lambda $ must look as shown in Figure~\ref{1curve} with
its single double point on the arc $D^{\prime }F$.

\begin{figure}[ht!]
\psfraga <-2pt,0pt> {D}{$D$}
\psfraga <-2pt,0pt> {E}{$E$}
\psfraga <-2pt,0pt> {F}{$F$}
\psfraga <-1pt,0pt> {H}{$H$}
\psfraga <-2pt,0pt> {G}{$G$}
\psfraga <-2pt,0pt> {g}{$\gamma$}
\centerline{\small
\includegraphics[width=.6\textwidth]{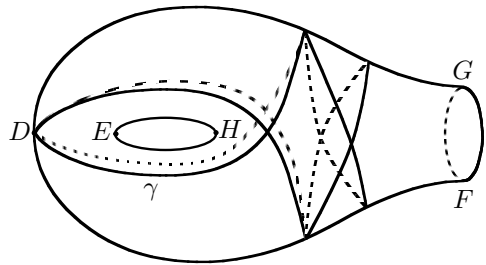}\hglue -2cm}
\caption{A curve with a unique configuration in any hyperbolic metric}
\label{1curve2}
\end{figure}

Now form a once punctured torus $T$ from $X$ by gluing together the two
boundary components containing $D,E$ and $D^{\prime },E^{\prime }$ so that $%
D $ is glued to $D^{\prime }$ and $E$ is glued to $E^{\prime }$. Then $\tau $
induces an orientation preserving involution on $T$ which we will continue
to denote by $\tau $, which fixes $D=D^{\prime },E=E^{\prime }$ and $H$. We
will be interested in the closed loop $\gamma $ on $T$ defined by $\gamma
=\lambda \cup \tau \lambda $. As $\tau $ is a rotation through $\pi $ in a
neighborhood of $D$, it follows that $\gamma $ is a closed geodesic. See
Figure~\ref{1curve2} which shows that $\gamma $ has seven double points. The
loop $\gamma $ has two innermost triangles, and using one of these triangles
we can change the configuration. However, we claim that no such triangle
move can be realized by the closed geodesics in a family of hyperbolic
metrics. For any hyperbolic metric on $X$ can be obtained from some all
right angled hyperbolic hexagon by doubling, so the preceding argument
applies to show that $\gamma $ will always have seven distinct double points
and no triple points. Thus the configuration of $\gamma $ cannot alter as
the hyperbolic metric changes continuously.

Next we discuss another example of a loop on a surface $F$ with several
minimal intersection configurations, of which only one is achieved by a
geodesic in a hyperbolic metric.

\begin{example}
A unique configuration on a thrice-punctured $S^2$.
\end{example}

\begin{figure}[ht!]
\centering
\begin{minipage}[c]{.49 \textwidth}
\centerline{\small
\psfraga <-2pt,0pt> {l}{$\lambda$}
\psfraga <-0pt,0pt> {v}{$v$}
\includegraphics[width=.7\textwidth]{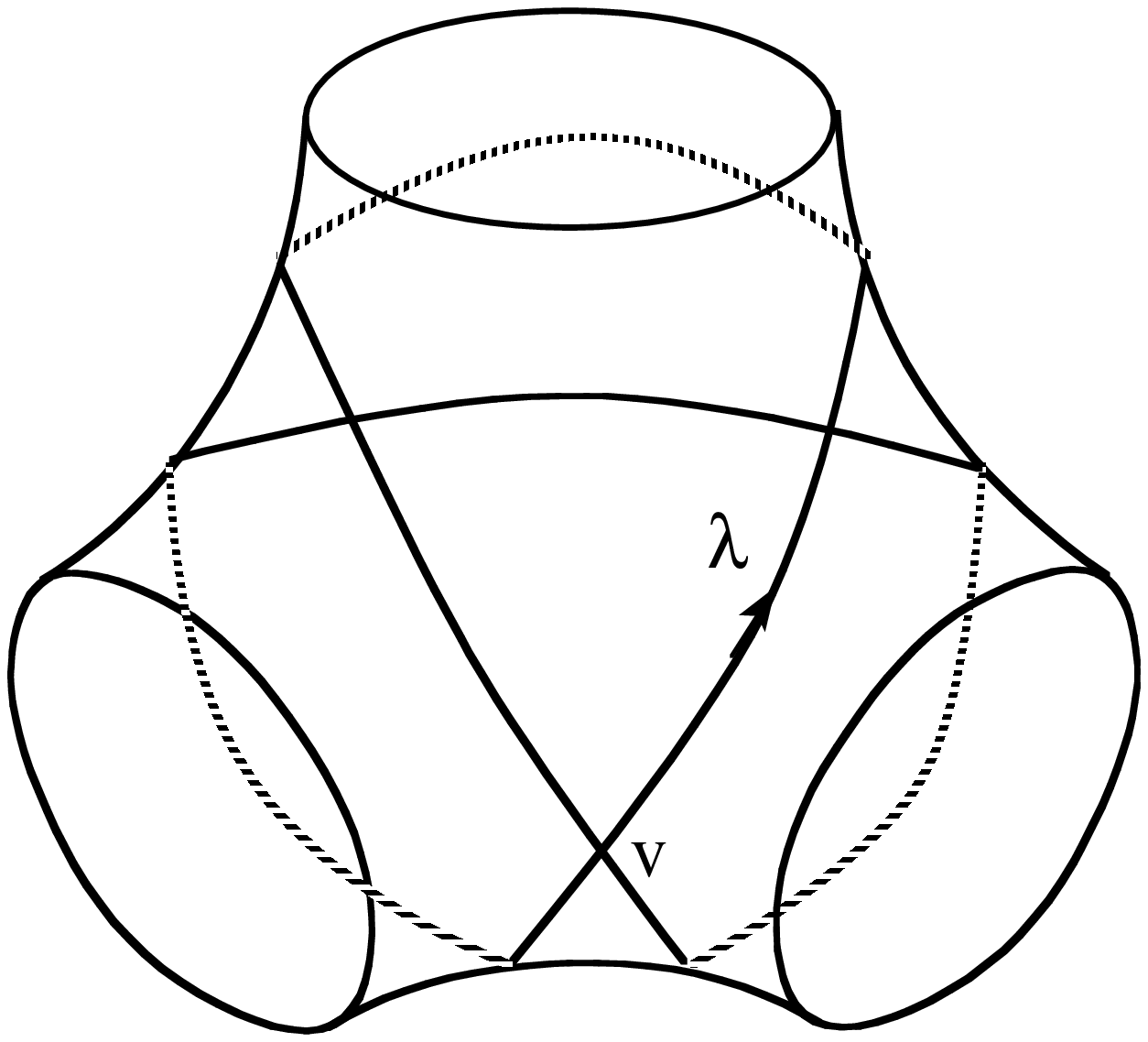}}
\end{minipage}
\begin{minipage}[c]{.49 \textwidth}
\centering
\includegraphics[width=.5\textwidth]{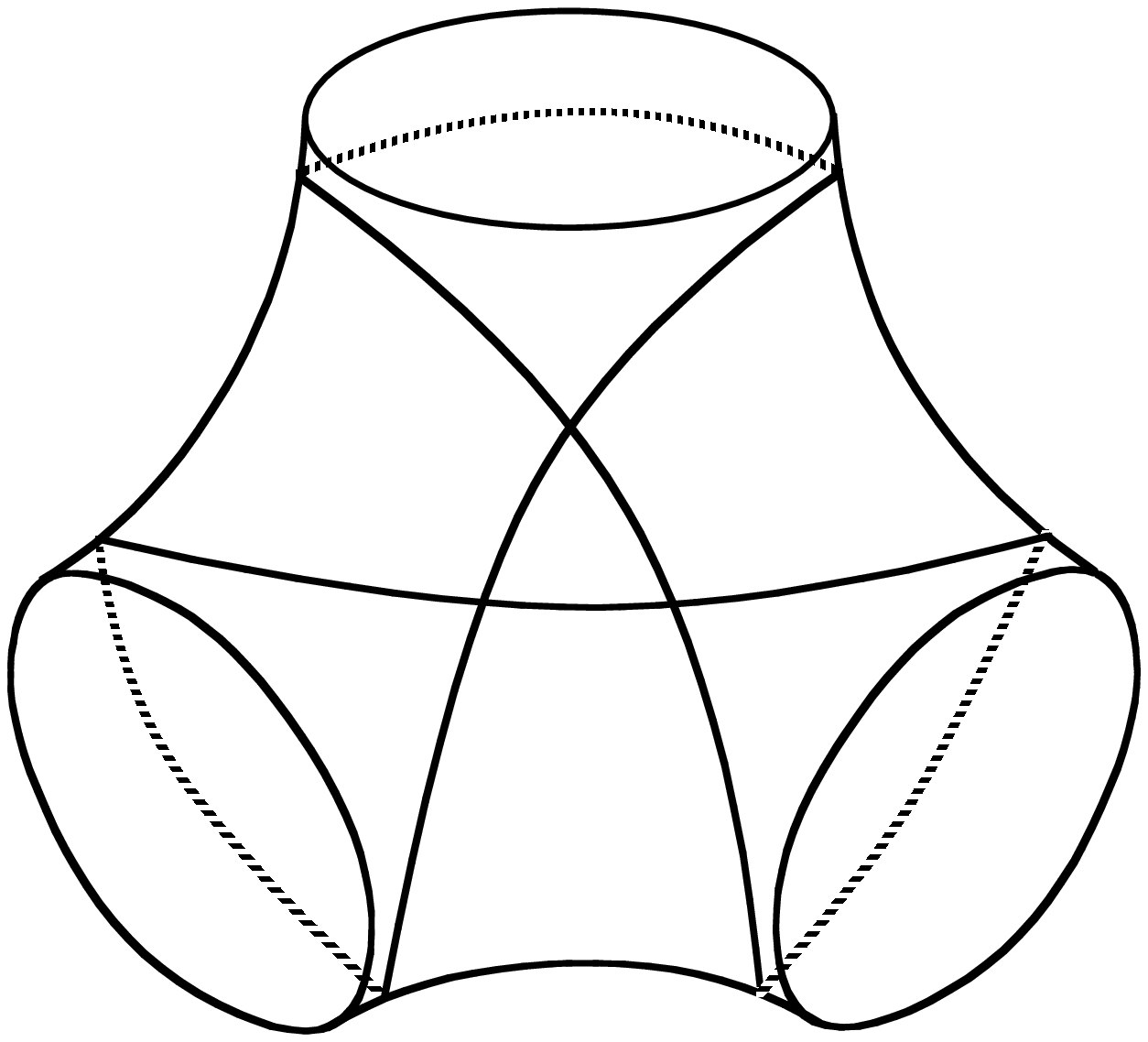}
\end{minipage}
\caption{A forced configuration and an impossible configuration}
\label{z3}
\end{figure}

Let $\Sigma $ denote a thrice-punctured $S^{2}$ equipped with a complete
hyperbolic metric of finite area. Let $\alpha $ denote the element of $\pi
_{1}(\Sigma )$ represented by the first loop shown in Figure~\ref{z3}. We
will use the fact that $\Sigma $ admits an action of $\mathbb{Z}_{3}$ by
isometries which cycles the three ends of $\Sigma $ to show that the
configuration of the closed geodesic representing $\alpha $ must be the
first one shown in Figure~\ref{z3}.

\begin{figure}[ht!]
\centering
\begin{minipage}[c]{.49 \textwidth}
\centerline{\small
\psfraga <-0pt,0pt> {l}{$\lambda$}
\psfraga <-2pt,-2pt> {v}{$v$}
\psfraga <-0pt,0pt> {q}{$\theta$}
\includegraphics[width=.7\textwidth]{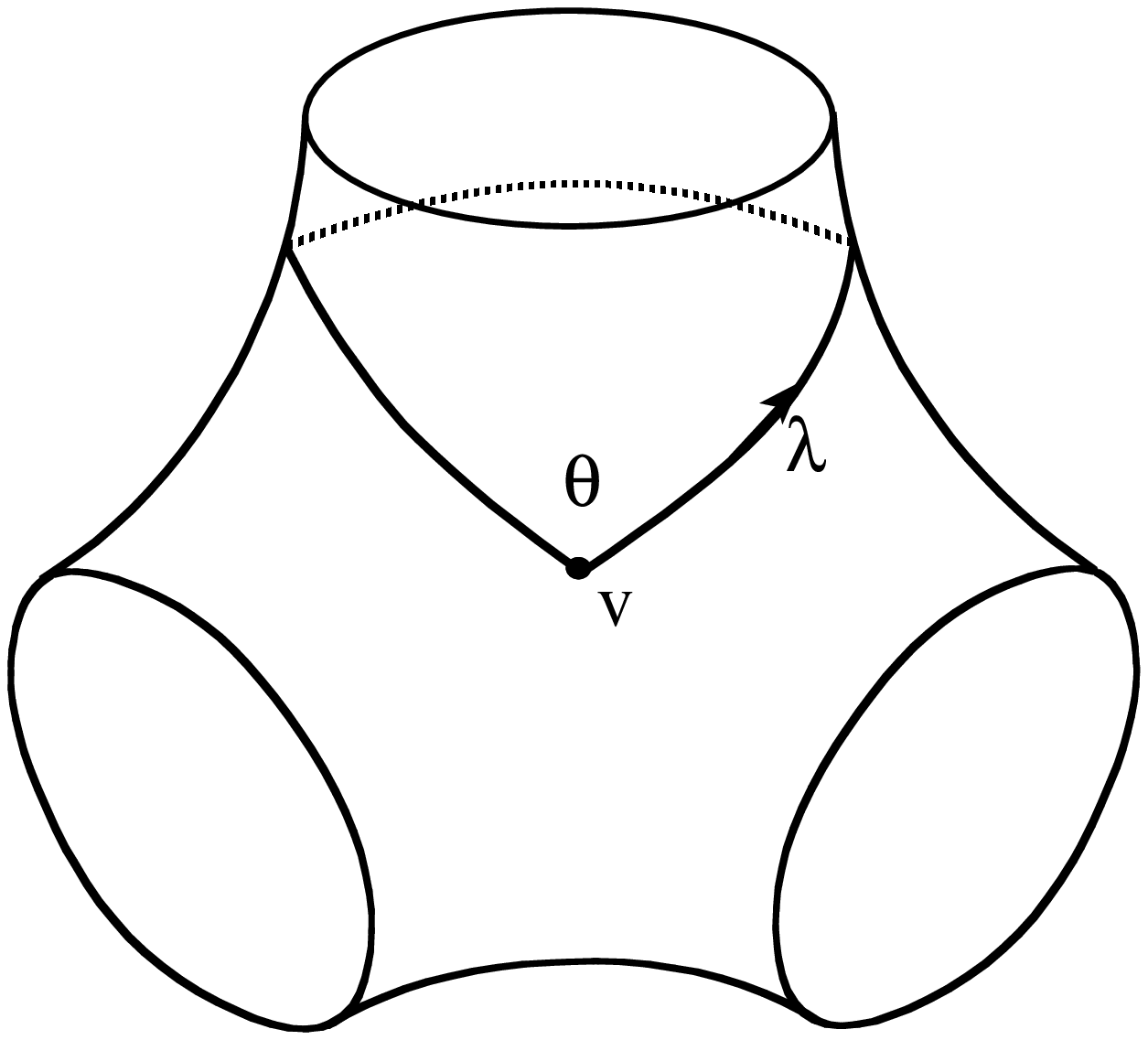}}
\end{minipage}
\begin{minipage}[c]{.49 \textwidth}
\centerline{
\includegraphics[width=.95\textwidth]{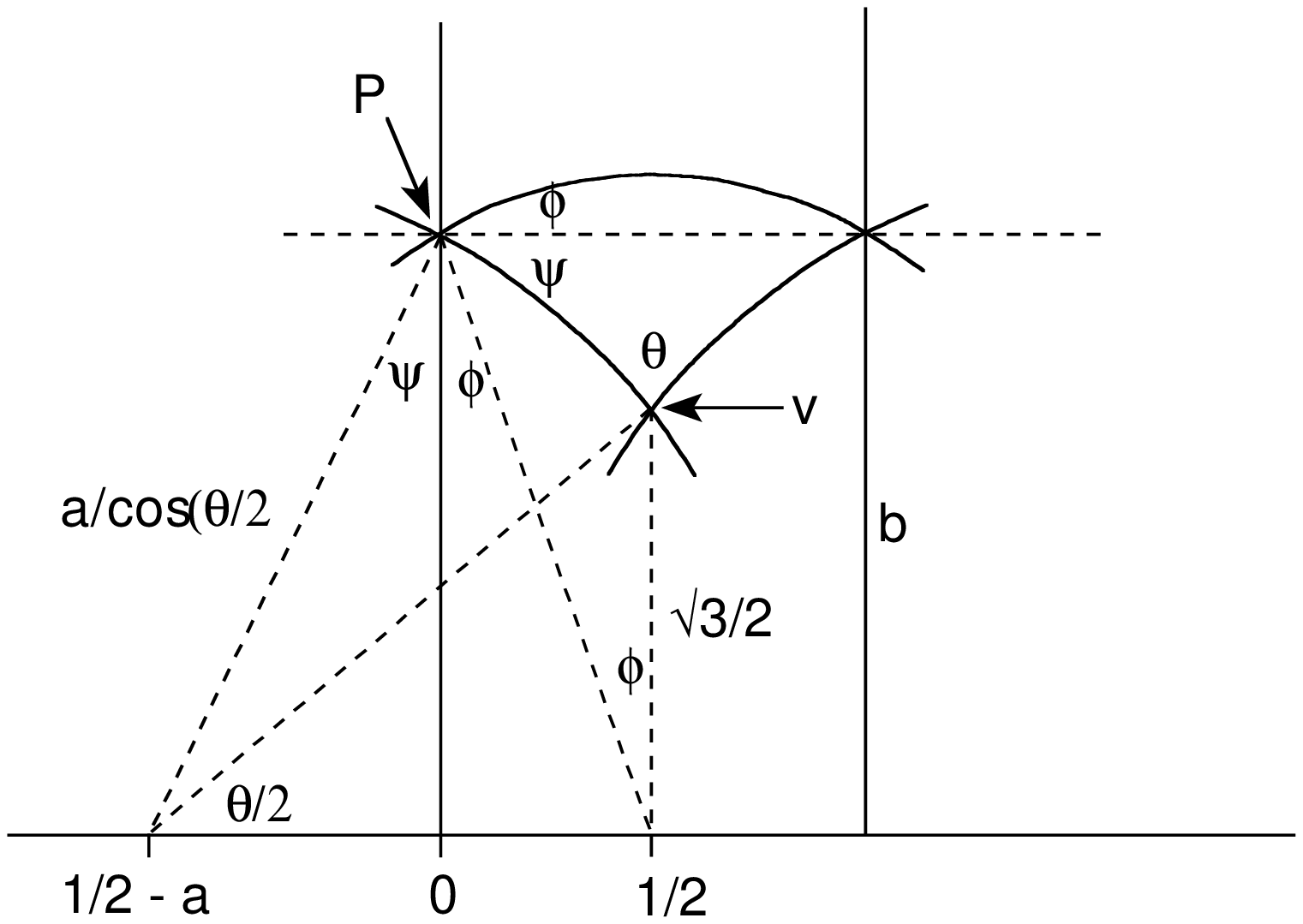}}
\end{minipage}
\caption{Calculating $\protect\theta$}
\label{theta}
\end{figure}

The proof is to consider the geodesic loop $\lambda $ shown in Figure~\ref
{theta}, whose corner is at $v$, one of the two points fixed by the action
of $\mathbb{Z}_{3}$, and show that $\theta >\pi /3$. Clearly the union of the
translates of $\lambda $ under the action of $\mathbb{Z}_{3}$ forms a loop
representing $\alpha $. If $\theta =\pi /3$, this loop will be a closed
geodesic and so the geodesic representing $\alpha $ will have a triple
point. If $\theta <\pi /3$, the geodesic representing $\alpha $ would have
the second configuration shown in Figure~\ref{z3}. The proof that $\theta
>\pi /3$ involves some straightforward hyperbolic geometry to show that $%
\theta =2\tan ^{-1}\left( \frac{\sqrt{3}}{2}\right) $, which is
approximately $81.79$ degrees. See Figure~\ref{theta}, which shows an ideal
triangle with vertices at $0,1$ and $\infty $ in the upper half plane model
of the hyperbolic plane. If we regard $\Sigma $ as the double of this
triangle, there is a natural quotient map from $\Sigma $ to the triangle and
the image of $\lambda $ is the piecewise geodesic triangular loop shown. It
has the properties that the exterior angles between $\lambda $ and the
geodesics $x=0$ and $x=1$ are all equal. Thus the angles marked $\phi $ and $%
\psi $ must be equal. We let $r$ denote the Euclidean radius of the circle
which forms the hyperbolic geodesic joining $v$ and $P$ and $a$ denote the
width of the projection to the $x$--axis of the radial segment of length $r$
connecting the center of this circle to $v$. Then $\displaystyle r =\frac{ a 
}{ \cos \theta /2}$ and the circle is centered at $( 1/2 - a, 0 )$. In these
coordinates, the rotation of the hyperbolic plane which sends $0$ to $1 $ to 
$\infty $ is the Mobius transformation $\displaystyle z\rightarrow \frac{1}{%
1-z}$. Recall that $v$ is fixed by this map. It follows that $v$ is the
point $\displaystyle\frac{1}{2}+i\frac{\sqrt{3}}{2}$. Hence $\displaystyle%
\tan \left( \theta /2\right) =\frac{\sqrt{3}}{2a}$. Also, if $b$ denotes the
Euclidean height of $P$ above the $x$--axis, then $\displaystyle\tan \phi =%
\frac{1/2}{b}$ and $\displaystyle\tan \psi =\frac{a-1/2}{b}$. As $\phi $ and 
$\psi $ are equal, we have $\displaystyle\tan \phi =\tan \psi , $ so that $%
a=1$. It follows that $\displaystyle\tan \left( \theta /2\right) =\frac{%
\sqrt{3}}{2}$, so that $\displaystyle\theta =2\tan ^{-1}\left( \frac{\sqrt{3}%
}{2}\right) $, as claimed.

\begin{rem}
Ian Agol has pointed out that the second configuration can also be
eliminated by a direct calculation in hyperbolic geometry. There is a
hexagon in the complement of the arcs, as well as a triangle. If $\theta
_{1},\theta _{2},\theta _{3}$ are the three interior angles of the triangle,
then the hexagon has exterior angles $\theta _{1},\theta _{2},\theta
_{3},\theta _{1},\theta _{2},\theta _{3}$. This is a contradiction since we
must have $\theta _{1}+\theta _{2}+\theta _{3}<\pi $ and $2\theta
_{1}+2\theta _{2}+2\theta _{3}>2\pi $. Moreover, Agol's observation applies
more generally in any complete negatively curved metric on the three
punctured sphere.
\end{rem}

Now we give an example of non-uniqueness of configurations.

\begin{example}
Non-unique configurations realized by a hyperbolic geodesic.
\end{example}

\begin{figure}[ht!]
\begin{minipage}{.49 \textwidth}
\centerline{\small 
\psfraga <-1pt, 0pt> {q}{$\theta$}
\psfraga <0pt, 0pt> {b}{$\beta$}
\includegraphics[width=1.5\textwidth]{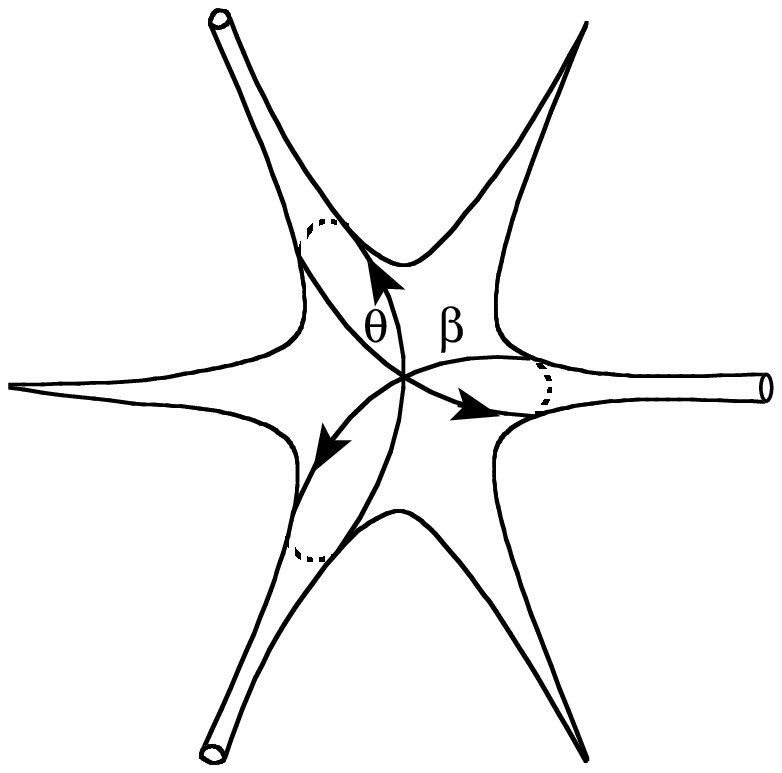}\hglue-5cm}
\end{minipage}
\begin{minipage}{.49 \textwidth}
\centerline{\small
\psfraga <0pt, -2pt> {q}{$\theta$}
\psfraga <-2pt, 0pt> {m}{$\mu$}
\includegraphics[width=1.2\textwidth]{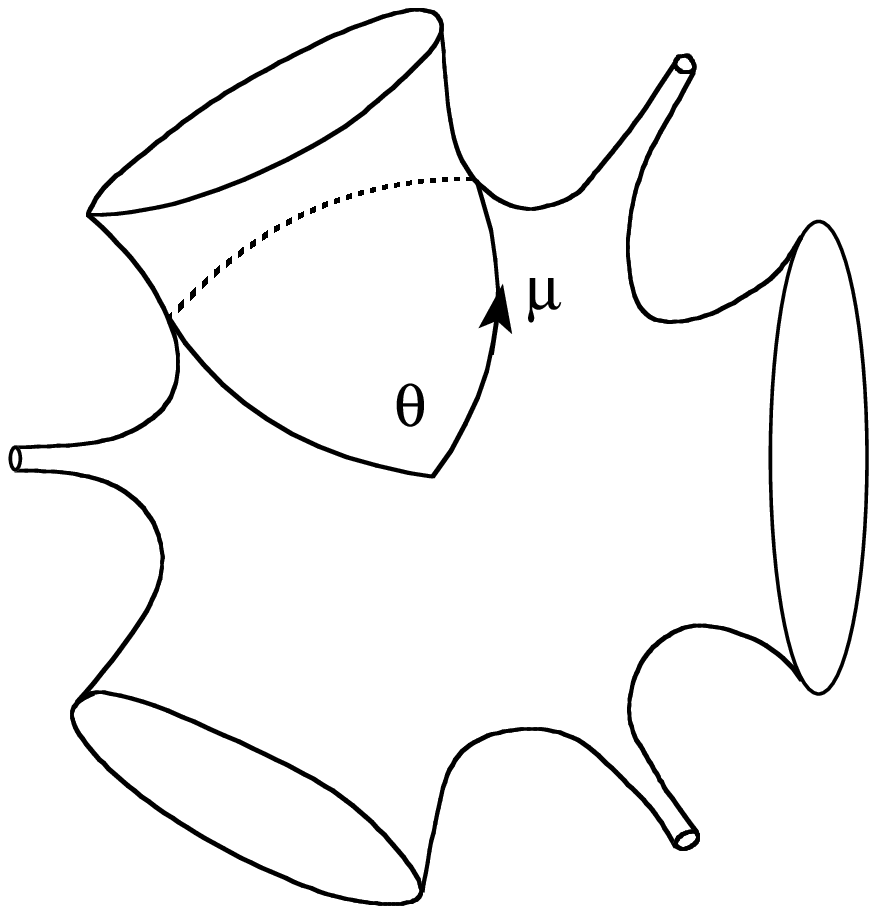}\hglue-3cm}
\end{minipage}
\caption{Two six punctured spheres with a hyperbolic metric}
\label{z6}
\end{figure}

Let $\Omega $ denote $S^{2}$ with six points removed equipped with a
complete hyperbolic metric having three cusp ends and three ends of infinite
area and admitting an action of $\mathbb{Z}_{3}$ which cycles the two types of end
among themselves. Let $\beta $ denote the element of $\pi _{1}(\Omega )$
represented by the loop shown in Figure~\ref{z6}. As before we consider the
arc $\mu $ shown in Figure~\ref{z6}, and the angle $\theta $. Note that the
union of the three translates of $\mu $ by the action of $\mathbb{Z}_{3}$ forms
a loop representing $\beta $. We will show that the closed geodesic
representing $\beta $ has at least two configurations which can be realized
by closed geodesics for some hyperbolic structure on $\Omega $.

To see this, start with a metric in which all the ends are cusps and there
is an action of $\mathbb{Z}_{6}$ on $\Omega $ which cycles the ends. In this
case, it is not as easy to calculate $\theta $, but it is clear that $\theta
<\pi /3$. Now alter the metric on $\Omega $, by enlarging the three
infinite area ends. Clearly $\theta \rightarrow 2\pi /3$ as the lengths of
the three closed geodesics tends to infinity. Hence by continuity, there is
a metric where $\theta =\pi /3$, and so the closed geodesic representing $%
\beta $ has a triple point. Distinct configurations will be obtained for
metrics near to this one for which $\theta <\pi /3$ and $\theta >\pi /3$.

Now we consider surfaces immersed in $3$--manifolds. There is a natural
analog of Theorem \ref{finite}. The statement, which we give below, uses the $1$--line
property, which is a property of a map of a surface into a $3$--manifold. See
sections 1 and 2 of \cite{HS2} for the definition and basic results.

\begin{thm}
\label{1-line} Let $N$ be a closed, $P^{2}$--irreducible $3$--manifold, let $F$
be a closed surface and let $f\co F\rightarrow N$ be a $2$--sided $\pi _{1}$%
--injective map. Then the general position immersions which are homotopic to $%
f$, have the 1--line property, whose double curves are primitive on $F$ and
have the least possible number of double points for their homotopy classes,
belong to only finitely many configurations.
\end{thm}

\begin{rem}
It is not assumed that $f$ has the 1--line property. Further, there may be no
immersions homotopic to $f$ with the required properties. In this case, the
result is trivial, but not interesting!
\end{rem}

\proof Homotopic maps with the 1--line property have precisely the
same double curves up to homotopy. Our hypothesis that the double curves
have the least number of double points means that we can use Theorem \ref{finite} to
deduce that there are only finitely many configurations for the double
curves of the maps homotopic to $f$ which have the 1--line property and the
other properties which we are assuming. Finally, the proof of Lemma 4.1 of 
\cite{HS2} shows that each configuration of double curves determines only
one configuration for a map $F^{2}\rightarrow N$, so the result
follows.\endproof

For any map $g$ of a surface into $N$, double points of the double curves of 
$g$ are triple points of $g$. Thus if $g$ satisfies all the hypotheses of
Theorem \ref{1-line}, then it must have the least possible number of triple
points in its homotopy class. However, the following example due to Casson
shows that a map may have the least possible number of triple points in its
homotopy class, while its double curves do not have the least possible
number of double points. In fact, {\em football regions}, complementary
regions homeomorphic to balls and bounded by three 2--gons, must occur in
Casson's example.

\begin{figure}[ht!]
\centerline{\small
\psfraga <-0pt,0pt> {C1}{$C_1$}
\psfraga <-0pt,0pt> {C2}{$C_2$}
\psfraga <-0pt,0pt> {C3}{$C_3$}
\psfraga <-2pt,0pt> {C4}{$C_4$}
\psfraga <-2pt,0pt> {C5}{$C_5$}
\psfraga <-2pt,0pt> {C6}{$C_6$}
\psfraga <-0pt,0pt> {45}{$\scriptstyle45$}
\psfraga <-1pt,0pt> {46}{$\scriptstyle46$}
\psfraga <-1pt,0pt> {56}{$\scriptstyle56$}
\includegraphics[width=.5\textwidth]{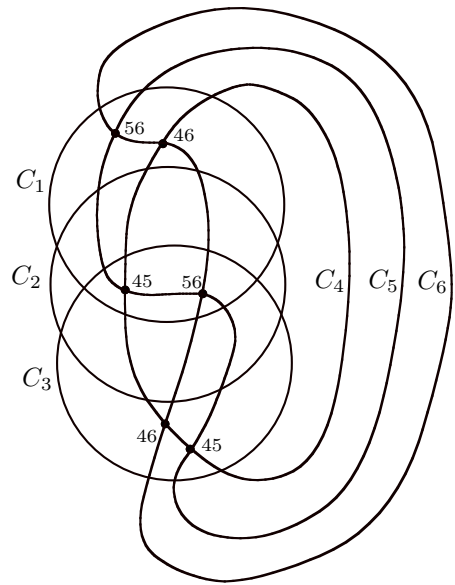}}
\caption{Any disks in a ball which have these curves as boundary must have a
football region between them.}
\label{casson}
\end{figure}

\begin{example}
(Casson)\qua A collection of surfaces which must contain a football region in
any configuration.
\end{example}

We consider six simple closed curves $C_{1},...,C_{6}$ on the 2--sphere $%
S^{2} $ as shown in Figure~\ref{casson}, so that each pair intersects
transversely in two points. Each $C_{i}$ bounds a 2--disc $D_{i}$ properly
embedded in the 3--ball $B^{3}$, and we assume that these discs are chosen in
general position. Further, by choosing these discs to be least area in some
metric, we can assume that any pair of these discs intersect in a single
arc, ie, there are no circles of intersection. The surprising property of
this picture is that there must be a football region $W$ in $B$, ie, a
sub-ball $W$ of $B$ bounded by the union of three discs each lying in some $%
D_{i}$, such that each pair of discs intersects in an arc. In particular, it
is impossible to embed the six discs $D_{i}$ in $B$ so that the double arcs
minimize their number of double points. Note that we are not claiming that $%
W $ is a component of the complement of the six discs. It is quite possible
that some of the discs can cut across $W$.

\begin{figure}[ht!]
\centering
\includegraphics[width=.5\textwidth]{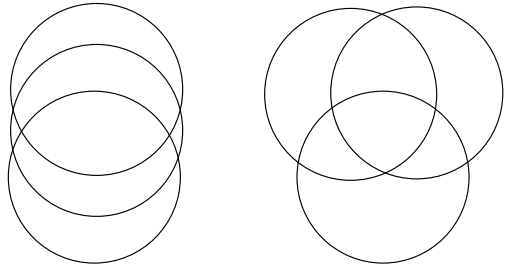}
\caption{Two configurations of three curves on a sphere}
\label{2config}
\end{figure}

Before starting on the proof, we remark that if one considers three simple
closed curves on $S^{2}$ which are in general position and such that each
pair intersect in exactly two points, then there are only two possible
configurations, as shown in Figure~\ref{2config}. In the first configuration
shown in Figure~\ref{2config}, which we refer to as the prism case, the
discs can be chosen so that each pair intersects in a single arc and there
is no triple point. In this case, the three discs cut $B^{3}$ into seven
regions, one of which meets $S^{2}$ in two triangular regions. This region
is referred to as the prism region. In the second configuration shown in
Figure~\ref{2config}, which we refer to as the triple point case, the discs
can be chosen so that each pair intersects in a single arc and there is
exactly one triple point. In the triple point case there must always be at
least one triple point however the discs are embedded. In Figure \ref{casson}
, the configuration of $C_{1},C_{2},C_{3}$ is of the prism type, and the
configuration of $C_{4},C_{5},C_{6}$ is of the triple point type.

Now suppose that we have an embedding of the $D_{i}$'s in $B^{3}$ such that
any two double lines of the $D_{i}$'s intersect in at most one point. Figure~%
\ref{casson} shows that the two ends of the double curve $D_{5}\cap D_{6}$
(labelled $56$ in the picture) lie on the same side of $D_{1}$, and that
this is on the opposite side of $D_{1}$ from the prism region $P$ formed by $%
D_{1},D_{2}$ and $D_{3}$. Similarly the two ends of the double curve $%
D_{4}\cap D_{5}$ (labelled $45$ in the picture) lie on the same side of $%
D_{3}$, and this is on the opposite side of $D_{3}$ from $P$. Finally the
two ends of the double curve $D_{4}\cap D_{6}$ (labelled $46$ in the
picture) lie on the same side of $D_{2}$, and this is on the opposite side
of $D_{2}$ from $P$. This implies that the three arcs in question cannot
have a common point, as the intersection of the sides of $D_{1},D_{2}$ and $%
D_{3}$ which do not contain the prism region $P$ is empty. This contradicts
the fact that the configuration of $D_{4},D_{5},D_{6}$ is of the triple
point type, so we conclude that for any embedding of the $D_{i}$'s in $B$
some pair of double lines $l$ and $m$ must intersect in at least two points.
For notational simplicity, suppose that $l=D_{1}\cap D_{2}$ and $m=D_{1}\cap
D_{3}$. Thus there are 2--gon regions in $D_{1}$ bounded by sub-arcs of $l$
and $m$. We choose one $X$ which is innermost in the sense that its interior
is disjoint from $l$ and $m$, and let $\lambda $ and $\mu $ denote the
sub-arcs of $l$ and $m$ respectively which form the boundary of $X$. Let $%
n=D_{2}\cap D_{3}$ and let $\nu $ denote the sub-arc of $n$ which has ends
at $\lambda \cap \mu $. Then $\lambda \cup \nu $ bounds a 2--gon $Y$ in $%
D_{2} $ and $\mu \cup \nu $ bounds a 2--gon $Z$ in $D_{3}$ and $X\cup Y\cup Z$
bounds a football region $W$ in $B^{3}$.

\rk{Acknowledgements}

Research of Joel Hass was partially supported by NSF grant DMS-9704286.
Research of Peter Scott was partially supported by NSF grant DMS-9626537.

\Addresses\recd

\end{document}